\documentclass[a4paper]{article}%
\usepackage{amssymb}
\usepackage{amsfonts}
\usepackage{amsmath}
\usepackage{geometry}
\usepackage{graphicx}%
\newtheorem{theorem}{Theorem}

\newtheorem{lemma}[theorem]{Lemma}

\newtheorem{remark}[theorem]{Remark}

\newenvironment{proof}[1][Proof]{\noindent\textbf{#1.} }{\ \rule{0.5em}{0.5em}}
\begin{document}

\title{Point Interaction Controls for the Energy Transfer in 3-D Quantum Systems.}
\author{Andrea Mantile\thanks{Dipartimento di Matematica e Applicazioni
\textquotedblleft R. Caccioppoli\textquotedblright\ Universit\`{a} degli Studi
Napoli \textquotedblleft Federico II\textquotedblright, Via Cintia, Monte S.
Angelo I-80126 Napoli, Italy. Supported through a European Community
Marie-Curie Fellowship and in the framework of CTS. E-mail:
andrea.mantile@dma.unina.it}}
\date{}
\maketitle

\begin{abstract}
We consider the problem of energy-mass transfer from scattering to bound
states for a one body quantum system subject to the action of a time dependent
point interaction in 3-D. Under suitable assumptions on the initial state of
the particle, we prove a result of local controllability of this process. Our
proof exploits the finite time asymptotic analysis of fractional integral
equations and the rank theorem for maps defined on Banach spaces.

\end{abstract}

\section{Introduction}

A perturbation of the Laplacian supported by a finite set of points $\left\{
y_{i}\right\}  _{i=1}^{n}$ in $\mathbb{R}^{d}$ - with $d\leq3$ - defines a
special case of singular perturbation referred to as point interaction. At a
formal level, the associated Schr\"{o}dinger operator can be written as:%
\begin{equation}
H=-\Delta+\sum_{i=1}^{n}\alpha_{i}\delta(x-y_{i}) \label{op_1}%
\end{equation}

These operators appeared first in Theoretical Physics during the 30's. They
were introduced in order to realize a model for the interaction of nucleons at
low energy \cite{fermi}. After, they became a natural tool to describe short
range forces or "small" obstacles for scattering of waves and particles.

From the point of view of applications, the main reason of interest of this
subject rests upon the fact that point interactions often lead to models which
are explicitly solvable. It turns out that the spectral characteristics
(eigenvalues and eigenfunctions) of operators (\ref{op_1}), and then all the
physical relevant quantities related to, can be explicitly computed
\cite{Albeverio1}. This circumstance motivates an increasing attention on the
application of point interaction models in various sciences, e.g. in physics,
chemistry, biology, and in technology.

In this work we build up a point interaction model of a time dependent
Schr\"{o}dinger operator; this interaction will be used as a control for the
energy transfer between continuous and discrete spectrum of a one body quantum
system. We will investigate the possibility of finding a time dependence
profile such that a prescribed part of the energy of a particle, initially
placed in a scattering state, moves on a bound state in finite time.

The subject we treat has its natural collocation in the framework of quantum
systems control theory. Our analysis may find applications in those contexts
where short range quantum potentials can be used as control tools.

\section{Point Interaction Model for Quantum Control Potentials and the Main
Result}

The rigorous definition of point interactions in Quantum Mechanics - due to
F.A. Berezin and L.D. Fadeev \cite{fad1} - rests upon the theory of
selfadjoint extensions of symmetric operators. The Hamiltonians describing
point interactions in the origin of $\mathbb{R}^{3}$ are defined by the
selfadjoint extensions of the symmetric operator%
\[
\left\{
\begin{array}
[c]{l}%
H=-\triangle\\
D(H)=C^{\infty}(\mathbb{R}^{3}\backslash\left\{  \underline{0}\right\}  )
\end{array}
\right.
\]
In \cite{Albeverio1} it has been shown that these extensions, denoted in the
following with $H_{\alpha}$, are parametrized by a real $\alpha\footnote{From
the physical point of view, the parameter $\alpha$ is linked to the inverse
scattering length of $H_{\alpha}$, \cite{Albeverio1}.}$; for any fixed
$\lambda\in\mathbb{R}^+$, the operator $H_\alpha$ can be represented
as\footnote{For the definitions and the properties concerning point
interactions operators, we refer to the book \cite{Albeverio1}}%
\begin{equation}
\left\{
\begin{array}
[c]{l}%
D(H_{\alpha})=\left\{  \psi\in L^{2}(\mathbb{R}^{3})\,\left\vert \,\psi
=\phi^{\lambda}+qG^{\lambda};\,\phi^{\lambda}\in H^{2}(\mathbb{R}^{3}%
);\,\phi^{\lambda}(\underline{0})=q(\alpha+\frac{\sqrt{\lambda}}{4\pi
})\right.  \right\} \\
\left(  H_{\alpha}+\lambda\right)  \psi=\left(  -\Delta+\lambda\right)
\phi^{\lambda}%
\end{array}
\right.  \label{Point Interaction 1}%
\end{equation}
where $G^{\lambda}$ is the Green function of $\left(  -\Delta+\lambda\right)
$%
\begin{equation}
G^{\lambda}(\underline{x})=\frac{e^{-\sqrt{\lambda}\,\left\vert \underline
{x}\right\vert }}{4\pi\,\left\vert \underline{x}\right\vert }
\label{Green Function}%
\end{equation}
As it follows from (\ref{Point Interaction 1}), every function in $D\left(
H_{\alpha}\right)  $ is composed of a regular part $\phi^{\lambda}\in H^{2}$
plus a singular term $G^{\lambda}$. It is worthwhile to notice that this
representation is not unique, but depends on the choice of $\lambda$. In
particular, for $\alpha<0$, we can fix: $\lambda=\left(  4\pi\left\vert
\alpha\right\vert \right)  ^{2}$, obtaining a null boundary condition:
$\phi^{\lambda}(\underline{0})=q(\alpha+\frac{\sqrt{\lambda}}{4\pi})=0$. In
this case the operator domain writes as%
\begin{equation}
D(H_{\alpha})=\left\{  \psi\in L^{2}(\mathbb{R}^{3})\left\vert \,\,\psi
=\phi+\frac{q}{4\pi\sqrt{2\left\vert \alpha\right\vert }}\psi_{\alpha}\right.
,\,\phi\in H^{2}(\mathbb{R}^{3}),\,\phi(\underline{0})=0,\,q\in\mathbb{C}%
\right\}  \label{dominio 1}%
\end{equation}
where $\psi_{\alpha}$ is the normalized Green function%
\begin{equation}
\psi_{\alpha}=\frac{\sqrt{2\left\vert \alpha\right\vert }e^{-4\pi\left\vert
\alpha\right\vert \,\left\vert \underline{x}\right\vert }}{^{\left\vert
\underline{x}\right\vert }} \label{bound state}%
\end{equation}
This representation, holding for negative values of the parameter $\alpha$,
will be extensively used in this work.

We recall the spectral properties of $H_{\alpha}$: the absolute continuous
spectrum coincides, for any $\alpha\in\mathbb{R}$, with the set $\left(
0,+\infty\right)  $; the point spectrum is empty for $\alpha\geq0$, or it
contains a single eigenvalue: $\lambda_{\alpha}=-16\pi^{2}\alpha^{2}$, for
$\alpha<0$. The related eigenstate is $\psi_{\alpha}$, defined in
(\ref{bound state}).

A time dependent point interaction Hamiltonian $H_{\alpha(t)}$ is defined by
fixing a real valued function, $\alpha(t)$, which characterizes the time
dependence profile. The Schr\"{o}dinger equation%
\begin{equation}
\left\{
\begin{array}
[c]{l}%
i\frac{d}{dt}\psi=H_{\alpha(t)}\psi\\
\psi(\underline{x},0)=\psi_{0}(\underline{x})\\
\psi(\underline{x},t)\in D(H_{\alpha(t)})
\end{array}
\right.  \label{Schroedinger_0}%
\end{equation}
describes the quantum dynamics generated by $H_{\alpha(t)}$. In \cite{Yafa1},
D.R. Yafaev has provided an explicit expression for the time propagator
associated to $H_{\alpha(t)}$ (equation (2.6) in the cited paper)%
\begin{equation}
\left\{
\begin{array}
[c]{l}%
\bigskip\psi(t,\underline{x})=U_{t}\psi_{0}(\underline{x})+i\int_{0}%
^{t}U(t-s,\underline{x})q(s)ds\\
q(t)+4\sqrt{\pi i}\int_{0}^{t}\frac{\alpha(s)q(s)}{\sqrt{t-s}}ds=4\sqrt{\pi
i}\int_{0}^{t}\frac{U_{s}\psi_{0}(\underline{0})}{\sqrt{t-s}}ds
\end{array}
\right.  \label{Schroedinger}%
\end{equation}
with%
\begin{equation}
U_{t}\psi(\underline{x})=\int_{\mathbb{R}^{3}}U(t,\underline{y}-\underline
{x})\psi(\underline{y})\,d\underline{y};\qquad U(t,\underline{x})=\frac
{1}{\left(  4\pi it\right)  ^{\frac{3}{2}}}e^{i\,\frac{\left\vert
\underline{x}\right\vert ^{2}}{4t}}\label{propagatore}%
\end{equation}
Given $\alpha\in L_{loc}^{\infty}(0,\infty)$, the map (\ref{Schroedinger}),
acting on the initial state $\psi_{0}$, is unitary in $L^{2}(\mathbb{R}^{3})$
and defines, in the weak sense, the solution of (\ref{Schroedinger_0}) at time
$t$ (see \cite{Yafa1} and references therein). The auxiliary variable $q(t)$,
solution of the second equation in (\ref{Schroedinger}), is usually referred
to as the \emph{charge} associated to the particle. It expresses the
coefficient of the singular part of the state $\psi(\underline{x},t)$. We
refer to the Appendix A for a study of the charge equation.

Let us fix, from now on, $\bar{\alpha}<0$. We assume, at $t=0$, the particle
to be placed in a scattering state of the Hamiltonian $H_{\bar{\alpha}}$%
\begin{equation}
\psi_{0}=\phi\in D(H_{\bar{\alpha}}):\left(  \phi,\psi_{\bar{\alpha}}\right)
_{L^{2}}=0 \label{1ip}%
\end{equation}
Our purpose is to find a suitable $\alpha(t)$ in the space of control
functions\footnote{The boundary conditions on $\alpha$ being chosen in order
to guarantee that at the inital and final times the Hamiltonian is
$H_{\bar{\alpha}}$.}%
\begin{equation}
\alpha\in H_{0}^{1}(0,T) \label{2ip}%
\end{equation}
such that the control interaction, $H_{\alpha(t)+\bar{\alpha}}$, is able to
steer the system (\ref{Schroedinger}) from $\psi_{0}$ to a state, $\psi(T)$,
whose projection $\left(  \psi(T),\psi_{\bar{\alpha}}\right)  $ has a
previously fixed value in a neighborhood of the origin of $\mathbb{C}$.

Replacing $\alpha$ with $\alpha+\bar{\alpha}$ in (\ref{Schroedinger}) and
projecting along $\psi_{\bar{\alpha}}$, we get a complex valued map%
\begin{equation}
\left\{
\begin{array}
[c]{l}%
\medskip F(\alpha):=\left(  U_{T}\phi+i\int_{0}^{T}U(T-s,\cdot)V(\alpha
)(s)ds,\,\psi_{\bar{\alpha}}\right)  _{L^{2}(\mathbb{R}^{3})}\\
V(\alpha)=q(t):q(t)+4\sqrt{\pi i}\int_{0}^{t}\frac{\left[  \alpha
(s)+\bar{\alpha}\right]  \,\,q(s)}{\sqrt{t-s}}ds=4\sqrt{\pi i}\int_{0}%
^{t}\frac{U_{s}\phi(\underline{0})}{\sqrt{t-s}}ds
\end{array}
\right.  \label{funzionale}%
\end{equation}
For $\alpha=0$, the system evolves under the action of $H_{\bar{\alpha}}$; in
this case we have%
\[
F(0)=\left(  e^{-iTH_{\bar{\alpha}}}\phi,\,\psi_{\bar{\alpha}}\right)
=e^{iT\lambda_{\bar{\alpha}}}\left(  \phi,\,\psi_{\bar{\alpha}}\right)  =0
\]
The condition for the solvability of our problem results to be the local
surjectivity of the map%
\begin{equation}
\left\{
\begin{array}
[c]{l}%
F(\alpha)=z\in\mathbb{C}\\
\alpha\in H_{0}^{1}(0,T)
\end{array}
\right.  \label{sistema}%
\end{equation}
around the point $\alpha=0$; i.e. we aim to prove that there exists a
neighborhood of the origin in $\mathbb{C}$, $I_{0}$, and a $C^{1}$-class map
$g:I_{0}\rightarrow H_{0}^{1}(0,T)$ such that%
\begin{equation}
F(g(z))=z\quad\forall z\in I_{0} \label{surjective}%
\end{equation}

The main goal of this work is to demonstrate the following result

\begin{theorem}
\label{Teorema 1}Let $\phi$ be a scattering state of the Hamiltonian
$H_{\bar{\alpha}}$ fulfilling the condition:%
\begin{equation}
\phi\in\left\{  \gamma\in S(\mathbb{R}^{3})\left\vert \,\left(  \gamma
,\psi_{\bar{\alpha}}\right)  _{L^{2}}=0;\,\gamma=\gamma(\left\vert
\underline{x}\right\vert );\,\gamma(0)=0\right.  \right\}
\label{stato iniziale}%
\end{equation}
being $S(\mathbb{R}^{3})$ the space of functions of rapid decrease (see e.g.
\cite{Simon}). Then the functional $F:H_{0}^{1}(0,T)\rightarrow\mathbb{C}$,
defined by (\ref{funzionale}), (\ref{2ip}), is a locally surjective map around
the point $\alpha=0$.
\end{theorem}

\begin{remark}
The inverse problem related to (\ref{funzionale}), i.e. find a control
$\alpha$ such that $F(\alpha)$ has a fixed value in $\mathbb{C}$, is connected
with the investigation of the energy exchanges in non autonomous quantum
systems. Although physical intuition suggests that the time dependence of the
Hamiltonian can force energy exchanges between the point and the continuous
spectrum, the dynamics describing these energy transfers is, in general,
rather complex (see, for instance, the study of the ionization problem driven
by time periodic delta interactions in \cite{Mik}-\cite{Cost1}); in
particular, we would like to stress that the local controllability of these
phenomena is not a simple consequence of the time dependence of the
Hamiltonian, but it depends from the specific structure of the interaction considered.
\end{remark}

We will adopt here a standard procedure in the analysis of nonlinear systems.
First we study the regularity properties of the map $F:H_{0}^{1}%
(0,T)\rightarrow\mathbb{C}$ and prove that $F\in C^{1}$(Lemma \ref{FC1}). Then
we investigate the surjectivity of the linearized map $d_{0}F$. To this aim we
will study a \emph{non controllability condition} for the linearized system
(Section 4); it will be shown that, under the hypothesis (\ref{stato iniziale}%
) on the initial state, this condition is never satisfied, obtaining, in this
way, a controllability result (sections 5-7). Then we conclude using a Rank
Theorem for functionals defined on Banach spaces (Section 8).

\section{The Nonlinear System}

In this section we investigate the differentiability properties of the
functional $F(\alpha)$. We shall denote with $X$ and $Y$ two Banach spaces,
$U$ an open subset of $X$, $d_{\alpha}F$ and $d_{\alpha}^{G}F$ respectively
the Fr\'{e}chet and the G\^{a}teaux derivatives of the map $F:U\rightarrow Y$
evaluated in the point $\alpha$. A differentiable functional $F:U\rightarrow
Y$ is said to be of class $C^{1}$ if the map:%
\begin{equation}
F^{\prime}:U\rightarrow L(X,Y),\ F^{\prime}(\alpha)=d_{\alpha}F \label{C1_def}%
\end{equation}
is continuous. Next we recall a standard result in the theory of differential
calculus in Banach spaces (Theorem 1.9 in \cite{Ambrosetti}):

\begin{theorem}
\label{prodi}Suppose $F:U\rightarrow Y$ is G\^{a}teaux-differentiable in $U$.
If the map:%
\[
F_{G}^{\prime}:U\rightarrow L(X,Y),\ F_{G}^{\prime}(\alpha)=d_{\alpha}^{G}F
\]
is continuous at $\alpha^{\ast}\in U$, then $F$ is Fr\'{e}chet-differentiable
at $\alpha^{\ast}$ and its Fr\'{e}chet derivative evaluated in $\alpha^{\ast}$
results:%
\[
d_{\alpha^{\ast}}F=d_{\alpha^{\ast}}^{G}F
\]

\end{theorem}

Consider the map $V:H_{0}^{1}(0,T)\rightarrow C(0,T)$ defined as follows%
\begin{equation}
\left[  V(\alpha)\right]  (t)=q(t):q(t)+4\sqrt{\pi i}\int_{0}^{t}\frac{\left[
\alpha(s)+\bar{\alpha}\right]  \,\,q(s)}{\sqrt{t-s}}ds=f(t)\label{map 0}%
\end{equation}
with $f\in C(0,T)$. From the estimate (\ref{stima_0}), in Lemma \ref{A}, we
have the following bound for the $C(0,T)$-norm of $\left[  V(\alpha)\right]
(t)$%
\begin{equation}
\left\Vert V(\alpha)\right\Vert _{C(0,T)}\leq\left\Vert f\right\Vert
_{C(0,T)}\Gamma(\alpha,T)\label{bound 1}%
\end{equation}
where $\Gamma(\alpha,T)$ is a positive finite constant depending on $\alpha$
and $T$.

\begin{lemma}
\label{VC1}Let $f\in C(0,T)$. The map $V:H_{0}^{1}(0,T)\rightarrow C(0,T)$
defined by (\ref{map 0}) is of class $C^{1}$.
\end{lemma}

\begin{proof}
First we prove that $V$ is continuous; let $\alpha,\beta\in H_{0}^{1}(0,T)$
and consider the difference $V(\alpha)(t)-V(\beta)(t)$: it solves the equation%
\[
q(t)+4\sqrt{\pi i}\int_{0}^{t}\frac{\left[  \alpha(s)+\bar{\alpha}\right]
\,\,q(s)}{\sqrt{t-s}}ds=4\sqrt{\pi i}\int_{0}^{t}\frac{\left[  \beta
(s)-\alpha(s)\right]  \,\left[  V(\beta)\right]  (s)}{\sqrt{t-s}}ds
\]
which is of type (\ref{map 0}); from the estimate (\ref{bound 1}) we obtain%
\begin{equation}
\left\Vert V(\alpha)-V(\beta)\right\Vert _{C(0,T)}\leq\left\vert 4\sqrt{\pi
i}\right\vert \left\Vert V(\beta)\right\Vert _{C(0,T)}\left\Vert \beta
-\alpha\right\Vert _{C(0,T)}2T^{\frac{1}{2}}\Gamma(\alpha,T) \label{stima1.1}%
\end{equation}
moreover, the Sobolev inequality%
\begin{equation}
\left\Vert g\right\Vert _{C(0,T)}\leq C\left\Vert g\right\Vert _{H^{1}(0,T)}
\label{stima1.2}%
\end{equation}
implies%
\begin{equation}
\left\Vert V(\alpha)-V(\beta)\right\Vert _{C(0,T)}\leq C\left\vert 4\sqrt{\pi
i}\right\vert \left\Vert V(\beta)\right\Vert _{C(0,T)}\left\Vert \beta
-\alpha\right\Vert _{H^{1}(0,T)}2T^{\frac{1}{2}}\Gamma(\alpha,T)
\label{stima1.3}%
\end{equation}
From (\ref{Gamma}), it follows that $\Gamma(\alpha,T)$ is a continuous
function: $C(0,T)\rightarrow\mathbb{R}^{+}$. Therefore $\Gamma(\alpha,T)$ is
uniformly bounded for $\alpha$ close to $\beta$ in $H^{1}$-norm. Previous
remark and the estimate (\ref{stima1.3}) imply the continuity of the
functional $V$ from $H_{0}^{1}(0,T)$ to $C(0,T)$.

Next we introduce the G\^{a}teaux derivative of $V$ in the point $\alpha$%
\[
d_{\alpha}^{G}V(u)=\lim_{\varepsilon\rightarrow0}\frac{V(\alpha+\varepsilon
u)-V\left(  \alpha\right)  }{\varepsilon}%
\]
this limit being in the sense of $C(0,T)$-topology. The explicit expression of
$d_{\alpha}^{G}V(u)$ is\footnote{This formula can be obtained by exploiting
the result of Lemma \ref{A}.}%
\begin{equation}
d_{\alpha}^{G}V(u)=q:q(t)+4\sqrt{\pi i}\int_{0}^{t}\frac{\left[
\alpha(s)+\bar{\alpha}\right]  \,\,q(s)}{\sqrt{t-s}}ds+4\sqrt{\pi i}\int
_{0}^{t}\frac{\left[  V(\alpha)\right]  (s)u(s)}{\sqrt{t-s}}ds=0 \label{dV}%
\end{equation}
This is a linear map $:H_{0}^{1}(0,T)\rightarrow C(0,T)$. Making use once more
of the estimate (\ref{stima_0}), we have%
\[
\left\Vert d_{\alpha}^{G}V(u)\right\Vert _{C(0,T)}\leq\left\vert 4\sqrt{\pi
i}\right\vert \left\Vert V(\alpha)\right\Vert _{C(0,T)}\left\Vert u\right\Vert
_{C(0,T)}\Gamma(\alpha,T)\,2T^{\frac{1}{2}}%
\]

In order to prove the statement of the Lemma we need to show that $V^{\prime
}:H_{0}^{1}(0,T)\rightarrow L(H_{0}^{1}(0,T),C(0,T))$, defined in
(\ref{C1_def}), is continuous. From Theorem (\ref{prodi}), it is sufficient to
prove that the map%
\[
V_{G}^{\prime}:H_{0}^{1}(0,T)\rightarrow L(H_{0}^{1}(0,T),C(0,T)),\ V_{G}%
^{\prime}(\alpha)=d_{\alpha}^{G}V
\]
is continuous, i.e.: $\lim_{\alpha\rightarrow\beta}\left\vert \left\vert
\left\vert d_{\alpha}^{G}V-d_{\beta}^{G}V\right\vert \right\vert \right\vert
=0$, where $\left\vert \left\vert \left\vert \cdot\right\vert \right\vert
\right\vert $ is the operator norm in $L(H_{0}^{1}(0,T),C(0,T))$. Set
$d_{\alpha}^{G}V(u)-d_{\beta}^{G}V(u)=q$; from the explicit formula
(\ref{dV}), it follows%
\[
q(t)+4\sqrt{\pi i}\int_{0}^{t}\frac{\left[  \alpha(s)+\bar{\alpha}\right]
\,\,q(s)}{\sqrt{t-s}}ds=4\sqrt{\pi i}\int_{0}^{t}\frac{\left\{  \left[
V(\beta)\right]  (s)-\left[  V(\alpha)\right]  (s)\right\}  u(s)-\left[
\alpha(s)-\beta(s)\right]  \left[  V(\beta)\right]  (s)}{\sqrt{t-s}}ds
\]
This equation is still of the type (\ref{map 0}); from the estimates
(\ref{bound 1}) and (\ref{stima1.2}) we get%
\[
\left\Vert q\right\Vert _{C(0,T)}\leq\left\vert 4\sqrt{\pi i}\right\vert
T^{\frac{1}{2}}C\left[  \left\Vert u\right\Vert _{H^{1}(0,T)}\left\Vert
V(\beta)-V(\alpha)\right\Vert _{C(0,T)}+\left\Vert \alpha-\beta\right\Vert
_{H^{1}(0,T)}\left\Vert f\right\Vert _{C(0,T)}\Gamma(\beta,T)\right]
\Gamma(\alpha,T)
\]
then the operator norm $\left\vert \left\vert \left\vert d_{\alpha}%
^{G}V-d_{\beta}^{G}V\right\vert \right\vert \right\vert $ can be bounded as
follows%
\begin{align*}
\left\vert \left\vert \left\vert d_{\alpha}^{G}V-d_{\beta}^{G}V\right\vert
\right\vert \right\vert  &  =\sup_{\substack{u\in H_{0}^{1}(0,T)\\\left\Vert
u\right\Vert =1}}\left\Vert d_{\alpha}^{G}V(u)-d_{\beta}^{G}V(u)\right\Vert
_{C(0,T)}\leq\\
&  \leq\left\vert 4\sqrt{\pi i}\right\vert T^{\frac{1}{2}}C\left[  \left\Vert
V(\beta)-V(\alpha)\right\Vert _{C(0,T)}+\left\Vert \alpha-\beta\right\Vert
_{H^{1}(0,T)}\left\Vert f\right\Vert _{C(0,T)}\Gamma(\beta,T)\right]
\Gamma(\alpha,T)
\end{align*}
As already noticed, $\Gamma(\alpha,T)$ is uniformly bounded when $\alpha$ is
close to $\beta$ in $H^{1}(0,T)$-norm. Then, the continuity of the map $V$
allows us to conclude that%
\[
\lim_{\left\Vert \alpha-\beta\right\Vert _{H_{0}^{1}(0,T)}\rightarrow
0}\left\vert \left\vert \left\vert d_{\alpha}^{G}V-d_{\beta}^{G}V\right\vert
\right\vert \right\vert =0
\]

\end{proof}

Our next task is the study of the regularity properties of the map defined in
(\ref{funzionale}) whose explicit expression is recalled here:%
\begin{equation}
\left\{
\begin{array}
[c]{l}%
\medskip F(\alpha):=\left(  U_{T}\phi(\underline{x})+i\int_{0}^{T}%
U(T-s,\underline{x})V(\alpha)(s)ds,\,\psi_{\bar{\alpha}}(\underline
{x})\right)  _{L^{2}(\mathbb{R}^{3})}\\
V(\alpha)=q(t):q(t)+4\sqrt{\pi i}\int_{0}^{t}\frac{\left[  \alpha
(s)+\bar{\alpha}\right]  \,\,q(s)}{\sqrt{t-s}}ds=4\sqrt{\pi i}\int_{0}%
^{t}\frac{U_{s}\phi(\underline{0})}{\sqrt{t-s}}ds
\end{array}
\right.  \label{map 1}%
\end{equation}

\begin{lemma}
\label{FC1}The map $F:H_{0}^{1}(0,T)\rightarrow\mathbb{C}$ defined by
(\ref{map 1}) is of class $C^{1}$.
\end{lemma}

\begin{proof}
First we show that $F$ is a continuous map. By inverting space and time
integrals in (\ref{map 1}), the scalar product $\left(  i\int_{0}%
^{T}U(T-s,\underline{x})V(\alpha)(s)ds,\,\psi_{\bar{\alpha}}(\underline
{x})\right)  _{L^{2}(\mathbb{R}^{3})}$ is given by%
\begin{align*}
\left(  i\int_{0}^{T}U(T-s,\underline{x})V(\alpha)(s)ds,\,\psi_{\bar{\alpha}%
}(\underline{x})\right)  _{L^{2}(\mathbb{R}^{3})}  &  =i\int_{0}^{T}%
\int_{\mathbb{R}^{3}}U(T-s,\underline{x})\psi_{\bar{\alpha}}(\underline
{x})d^{3}\underline{x}\,V(\alpha)(s)ds=\\
&  =i\int_{0}^{T}U_{T-s}\psi_{\bar{\alpha}}(\underline{0})\,V(\alpha)(s)ds
\end{align*}
Let $\alpha,\beta\in H_{0}^{1}(0,T)$; making use of the previous expression,
the difference $F(\alpha)-F(\beta)$ can be written as%
\[
F(\alpha)-F(\beta)=i\int_{0}^{T}U_{T-s}\psi_{\bar{\alpha}}(\underline
{0})\,\left(  V(\alpha)(s)-V(\beta)(s)\right)  \,ds
\]
Therefore, it satisfies the estimate%
\begin{equation}
\left\vert F(\alpha)-F(\beta)\right\vert \leq\left\Vert V(\alpha
)(s)-V(\beta)(s)\right\Vert _{C(0,T)}\int_{0}^{T}\left\vert U_{T-s}\psi
_{\bar{\alpha}}(\underline{0})\right\vert \,ds \label{stima4}%
\end{equation}
The continuity of $F$ then follows directly from Lemma \ref{VC1}. Moreover,
from Theorem \ref{prodi} and definition (\ref{C1_def}) it follows that $F$ is
of $C^{1}$ class if the map%
\[
F_{G}^{\prime}:H_{0}^{1}(0,T)\rightarrow L(H_{0}^{1}(0,T),\mathbb{C}%
),\ F_{G}^{\prime}(\alpha)=d_{\alpha}^{G}F
\]
is continuous, i.e.: $\lim_{\alpha\rightarrow\beta}\left\vert \left\vert
\left\vert d_{\alpha}^{G}F-d_{\beta}^{G}F\right\vert \right\vert \right\vert
=0$. The G\^{a}teaux derivative of $F$ evaluated in the point $\alpha$ and
acting on $u$ is%
\begin{equation}
d_{\alpha}^{G}F(u)=i\left(  \int_{0}^{T}U(T-s,\underline{x})\,d_{\alpha
}V(u)(s)ds,\,\psi_{\bar{\alpha}}(\underline{x})\right)  _{L^{2}(\mathbb{R}%
^{3})} \label{dF}%
\end{equation}
By inverting space and time integrals, we obtain%
\begin{equation}
d_{\alpha}^{G}F(u)=i\int_{0}^{T}U_{T-s}\psi_{\bar{\alpha}}(\underline
{0})\,d_{\alpha}^{G}V(u)(s)ds \label{dF 1}%
\end{equation}
Consider the difference $d_{\alpha}^{G}F(u)-d_{\beta}^{G}F(u)$; from
(\ref{dF 1}) it results
\begin{align*}
\left\vert d_{\alpha}^{G}F(u)-d_{\beta}^{G}F(u)\right\vert  &  =\left\vert
\int_{0}^{T}U_{T-s}\psi_{\bar{\alpha}}(\underline{0})\,\left(  d_{\alpha}%
^{G}V(u)(s)-d_{\beta}^{G}V(u)(s)\right)  ds\right\vert \leq\\
&  \leq\left\Vert d_{\alpha}^{G}V(u)-d_{\beta}^{G}V(u)\right\Vert
_{C(0,T)}\int_{0}^{T}\left\vert U_{T-s}\psi_{\bar{\alpha}}(\underline
{0})\right\vert \,ds
\end{align*}
and%
\[
\left\vert \left\vert \left\vert d_{\alpha}^{G}F-d_{\beta}^{G}F\right\vert
\right\vert \right\vert =\sup_{\substack{u\in H_{0}^{1}(0,T)\\\left\Vert
u\right\Vert =1}}\left\vert d_{\alpha}^{G}F(u)-d_{\beta}^{G}F(u)\right\vert
\leq\left\vert \left\vert \left\vert d_{\alpha}^{G}V-d_{\beta}^{G}V\right\vert
\right\vert \right\vert \int_{0}^{T}\left\vert U_{T-s}\psi_{\bar{\alpha}%
}(\underline{0})\right\vert \,ds
\]
The continuity of the map $F_{G}^{\prime}$ easily follows from the continuity
of the map $V_{G}^{\prime}$.
\end{proof}

\section{The Linearized System}

We consider the map:%
\begin{equation}
\left\{
\begin{array}
[c]{l}%
d_{0}F(u)=z\\
u\in H_{0}^{1}(0,T)
\end{array}
\right.  \label{d-sistema}%
\end{equation}
where $d_{0}F$ is the Fr\'{e}chet derivative of functional $F$ evaluated in
$\alpha=0$. Our aim is to prove the following result:

\begin{theorem}
\label{Teorema 1.1}Under the assumptions of Theorem \ref{Teorema 1}, the map
defined by (\ref{d-sistema}) is surjective.
\end{theorem}

The proof of Theorem \ref{Teorema 1.1} will be given in Sections 4-7 following
an \emph{ad absurdum }procedure.

First we write the functional (\ref{d-sistema}) into an explicit form:%

\begin{equation}
\left\{
\begin{array}
[c]{l}%
\bigskip d_{0}F(u)=i\left(  \int_{0}^{T}U(T-s,\cdot)\,d_{0}V(u)(s)ds,\,\psi
_{\bar{\alpha}}\right)  _{L^{2}(\mathbb{R}^{3})}\\
\,d_{0}V(u)=q(t):q(t)+4\sqrt{\pi i}\,\bar{\alpha}\int_{0}^{t}\frac
{\,q(s)}{\sqrt{t-s}}ds=-4\sqrt{\pi i}\int_{0}^{t}\frac{u(s)\,V(0)(s)}%
{\sqrt{t-s}}ds
\end{array}
\right.  \label{dF0}%
\end{equation}
Here the dependence of $d_{0}F(u)$ from $V(0)(t)$ may be emphasized by making
use of the relation (\ref{soluzione 2}) (in Appendix A), plus the Fubini
Theorem (in order to exchange time and space integrations in (\ref{dF0}));
proceeding in this way, we get%
\begin{multline*}
d_{0}F(u)=i\int_{0}^{T}U_{T-s}\psi_{\bar{\alpha}}(\underline{0})\,d_{0}%
V(u)(s)ds=\\
=4\pi i^{\frac{3}{2}}\int_{0}^{T}U_{T-s}\psi_{\bar{\alpha}}(\underline
{0})\,\int_{0}^{s}G(t-s^{\prime})\,V(0)(s^{\prime})\,u(s^{\prime}%
)\,ds^{\prime}ds
\end{multline*}
Applying the Dirichlet formula to the double integral we finally obtain%
\begin{equation}
d_{0}F(u)=4\pi i^{\frac{3}{2}}\int_{0}^{T}\,V(0)(s^{\prime})\,u(s^{\prime
})\int_{s^{\prime}}^{T}U_{T-s}\psi_{\bar{\alpha}}(\underline{0}%
)\,G(s-s^{\prime})ds\,ds^{\prime} \label{dF0(1)}%
\end{equation}

Let us suppose the map $d_{0}F:H_{0}^{1}(0,T)\rightarrow\mathbb{C}$ to be
non-surjective; then $d_{0}F(u)$ should have a constant direction in the
complex plane for any $u\in H_{0}^{1}(0,T)$. We can express this as a
\textbf{non-controllability condition}%
\[
\exists\;C\in\mathbb{C}:\left\{
\begin{array}
[c]{l}%
C\neq0\\
C\cdot\int_{0}^{T}V(0)(s^{\prime})\,u(s^{\prime})\int_{s^{\prime}}^{T}%
U_{T-s}\psi_{\bar{\alpha}}(\underline{0})\,G(s-s^{\prime})ds\,ds^{\prime}=0
\end{array}
\right.  \quad\forall u\in H_{0}^{1}(0,T)
\]
where '$\cdot$' indicates the scalar product in $\mathbb{C}$. In particular,
being this condition true for any real valued $u\in C_{0}^{\infty}(0,T)$, the
fundamental lemma of calculus of variations implies%
\begin{equation}
\exists\;C\in\mathbb{C}:\left\{
\begin{array}
[c]{l}%
C\neq0\\
C\cdot V(0)(t)\int_{t}^{T}U_{T-s}\psi_{\bar{\alpha}}(\underline{0}%
)\,G(s-t)ds=0
\end{array}
\right.  \quad\forall t\in\left[  0,T\right]  \label{condizione 0}%
\end{equation}
Making use of (\ref{carica inversa 2}) - in Appendix A - we get%
\begin{equation}
\exists\;C\in\mathbb{C}:\left\{
\begin{array}
[c]{l}%
C\neq0\\
C\cdot V(0)(t)e^{-i(T-t)\lambda_{\bar{\alpha}}}=0
\end{array}
\right.  \quad\forall t\in\left[  0,T\right]  \label{condizione 1}%
\end{equation}

\section{Small Time Asymptotic for the Charge}

In this Section we consider the asymptotic behavior of the charge $V(0)$ for
$t\rightarrow0$. We shall denote with $o(t)$ and $O(t)$ two complex valued
functions of the real variable $t$ satisfying the conditions%
\begin{equation}
\lim_{t\rightarrow0}\frac{o(t)}{t}=0
\end{equation}%
\begin{equation}
\lim\sup_{t\rightarrow0}\left\vert \frac{O(t)}{t}\right\vert <\infty\label{O}%
\end{equation}

According to our hypothesis on the initial state, (\ref{stato iniziale}), and
to definition (\ref{funzionale}), the function $V(0)$ is the solution of the
equation%
\[
q(t)+4\sqrt{\pi i}\bar{\alpha}\int_{0}^{t}\frac{\,q(s)}{\sqrt{t-s}}%
ds=4\sqrt{\pi i}\int_{0}^{t}\frac{U_{s}\gamma(\underline{0})}{\sqrt{t-s}}ds
\]
therefore, thanks to (\ref{soluzione 1}), it may be represented in the form%
\begin{equation}
V(0)(t)=4\pi\sqrt{i}\int_{0}^{t}G(t-s)U_{s}\gamma(\underline{0})\,ds
\label{carica}%
\end{equation}
Its small time behavior is connected to the limiting behavior of $G(t)$ and
$U_{t}\gamma(\underline{0})$ for $t\rightarrow0$. In order to study this
problem we need the following Lemmas

\begin{lemma}
\label{lemma 1}Let $\gamma=\gamma(\left\vert \underline{x}\right\vert )$ be a
radial function belonging to the space of functions of rapid decrease
$S(\mathbb{R}^{3})$ (see e.g. \cite{Simon}). If we assume that%
\begin{equation}
D_{0}=\left\{  n\in\mathbb{N}:\left.  \frac{d^{n}}{dt^{n}}U_{t}\gamma
(\underline{0})\right\vert _{t=0}\neq0\right\}  \label{3ip}%
\end{equation}
is a non empty set, then the function $U_{t}\gamma(\underline{0})$ admits the
power expansion:%
\begin{equation}
U_{t}\gamma(\underline{0})=a_{m}\,t^{m}+O(t^{m+1});\ a_{m}\neq0
\label{sviluppo 1}%
\end{equation}%
\begin{equation}
a_{m}=\frac{4\pi}{\left(  2\pi\right)  ^{\frac{3}{2}}}\frac{\left(  -i\right)
^{m}}{m!}\int_{0}^{+\infty}k^{2m+2}\,\mathcal{F}\gamma(k)\,dk
\label{sviluppo 1.1}%
\end{equation}
where $\mathcal{F}$ denotes the Fourier transform, $m=\min D_{0}$ and $O(t)\in
C^{\infty}[0,+\infty)$.
\end{lemma}

\begin{proof}
The Fourier transform operator, $\mathcal{F}$, is an homeomorphism of the
space $S(\mathbb{R}^{3})$ in itself. It acts on $U_{t}\gamma(\underline{0})$
as follows:%
\[
\mathcal{F}U_{t}\gamma(k)=\mathcal{F}\gamma(k)\,e^{-ik^{2}t}%
\]
Then, using $\mathcal{F}^{-1}$, we can represent $U_{t}\gamma(\underline{0})$
in the form%
\begin{equation}
U_{t}\gamma(\underline{0})=\frac{4\pi}{\left(  2\pi\right)  ^{\frac{3}{2}}%
}\int_{0}^{+\infty}k^{2}\mathcal{F}\gamma(k)\,e^{-ik^{2}t}\,dk\label{lem1-1}%
\end{equation}
From the regularity assumptions on $\gamma$, we have: $\mathcal{F}\gamma(k)\in
S(\mathbb{R}^{3})$ and $U_{t}\gamma(\underline{0})\in C^{\infty}[0,+\infty)$.
Setting $m=\min D_{0}$, the Taylor's expansion of $U_{t}\gamma(\underline{0})$
up to order $m$ in a right neighborhood of the origin is%
\[
U_{t}\gamma(\underline{0})=a_{m}\,t^{m}+O(t^{m+1})
\]
with $O(t)\in C^{\infty}[0,+\infty)$. The coefficient $a_{m}$ is different
from zero - due to the hypothesis (\ref{3ip}) - and explicitly given by
$\left.  \frac{1}{m!}\frac{d^{m}}{dt^{m}}U_{t}\gamma(\underline{0})\right\vert
_{t=0}$; form (\ref{lem1-1}), this quantity is%
\[
a_{m}=\frac{4\pi}{\left(  2\pi\right)  ^{\frac{3}{2}}}\frac{\left(  -i\right)
^{m}}{m!}\int_{0}^{+\infty}k^{2m+2}\,\mathcal{F}\gamma(k)\,dk
\]

\end{proof}

\begin{lemma}
\label{lemma 2}Let $G(t)$ be given by (\ref{kernel});\ then it admits the
following representation%
\begin{equation}
G(t)=\frac{1}{\sqrt{\pi t}}+b_{0}+O(t);\quad b_{0}=-4\pi\bar{\alpha}\sqrt{i}
\label{sviluppo 2}%
\end{equation}
where $O(t)\in C^{\infty}[0,+\infty)$.
\end{lemma}

\begin{proof}
The proof easily follows from the analytic properties of the $erfc(t)$
function (\cite{Abramowitz}, relation 7.1.5, pag 297 )
\end{proof}

We will use these results to get an expansion in power of $t^{\frac{1}{2}}$
for the charge (\ref{carica}). If we assume condition (\ref{3ip}) to hold,
Lemma \ref{lemma 1} may be applied to our case, with the only restriction
$m\neq0$ due to the fact that the boundary condition $\gamma(\underline{0})=0$
implies $U_{t}\gamma(\underline{0})=0$ for $t=0$. By substitution of
(\ref{sviluppo 1}) and (\ref{sviluppo 2}) into (\ref{carica}) we obtain%
\begin{multline*}
V(0)(t)=4\sqrt{\pi i}\,a_{m}\int_{0}^{t}\frac{s^{m}}{\sqrt{t-s}}ds+4\pi
\sqrt{i}\,a_{m}b_{0}\int_{0}^{t}s^{m}ds+4\sqrt{\pi i}\int_{0}^{t}%
\frac{O(s^{m+1})}{\sqrt{t-s}}ds+\\
+4\pi\sqrt{i}\,b_{0}\int_{0}^{t}O(s^{m+1})ds+4\pi\sqrt{i}\,a_{m}\int_{0}%
^{t}O(t-s)\,s^{m}ds+4\pi\sqrt{i}\int_{0}^{t}O(t-s)\,O(s^{m+1})ds
\end{multline*}
An explicit calculation of the first terms in this expression leads us to the
following expansion%
\begin{multline*}
V(0)(t)=A_{m}a_{m}\,\sqrt{i}\,t^{m+\frac{1}{2}}+B_{m}a_{m}\,\sqrt{i}%
\,b_{0}\,t^{m+1}+4\sqrt{\pi i}\int_{0}^{t}\frac{O(s^{m+1})}{\sqrt{t-s}}ds+\\
+4\pi\sqrt{i}\,b_{0}\int_{0}^{t}O(s^{m+1})ds+4\pi\sqrt{i}\,a_{m}\int_{0}%
^{t}O(t-s)\,s^{m}ds+4\pi\sqrt{i}\int_{0}^{t}O(t-s)\,O(s^{m+1})ds
\end{multline*}
where $A_{m}$ and $B_{m}$ are strictly positive real constants. By definition
(\ref{O}), it exists $C$ such that $\left\vert O(s^{k})\right\vert \leq
C\left\vert s^{k}\right\vert \leq C\left\vert t^{k}\right\vert $ for any
$\left\vert s\right\vert \leq\left\vert t\right\vert $; using this estimate,
the following relations are easily obtained%
\begin{equation}
\int_{0}^{t}O(s^{m+1})ds=O(t^{m+2})\label{1}%
\end{equation}%
\begin{equation}
\int_{0}^{t}O(t-s)\,s^{m}ds=O(t^{m+2})\label{2}%
\end{equation}%
\begin{equation}
\int_{0}^{t}O(t-s)\,O(s^{m+1})ds=o(t^{m+2})\label{3}%
\end{equation}%
\begin{equation}
\int_{0}^{t}\frac{O(s^{m+1})}{\sqrt{t-s}}ds=O(t^{m+\frac{3}{2}})\label{4}%
\end{equation}
and the small time asymptotic representation for $V(0)$ can be rewritten as%
\begin{equation}
V(0)(t)=A_{m}a_{m}\,\sqrt{i}\,t^{m+\frac{1}{2}}+B_{m}a_{m}\,\sqrt{i}%
\,b_{0}\,t^{m+1}+O(t^{m+\frac{3}{2}})\label{sviluppo 3}%
\end{equation}

\section{The Non Controllability Condition in the limit $t\rightarrow0$}

Here we study condition (\ref{condizione 1}) in a neighborhood $\left[
0,\delta\right)  $ of the origin with $\delta<T$.

Let us first set the condition (\ref{condizione 1}) in the equivalent form%
\begin{equation}
\exists\,K\in\left[  0,2\pi\right)  :\arg\left(  V(0)(t)e^{-i(T-t)\lambda
_{\bar{\alpha}}}\right)  =K\quad\forall t\in\lbrack0,T] \label{condizione 1.1}%
\end{equation}
where $\arg z$ is defined modulus $2\pi$. Making use of (\ref{carica}), we see
that%
\begin{gather*}
\arg\left(  V(0)(t)e^{-i(T-t)\lambda_{\bar{\alpha}}}\right)  =\qquad
\qquad\qquad\qquad\qquad\qquad\qquad\qquad\qquad\\
=\arg\left(  4\pi\sqrt{i}\int_{0}^{t}G(t-s)U_{s}\gamma(\underline
{0})\,ds\,e^{-i(T-t)\lambda_{\bar{\alpha}}}\right)  \qquad
\end{gather*}
Then, condition (\ref{condizione 1.1}) implies%
\begin{equation}
\exists\,K\in\left[  0,2\pi\right)  :\arg\left(  4\pi\sqrt{i}\int_{0}%
^{t}G(t-s)U_{s}\gamma(\underline{0})\,ds\,e^{-i(T-t)\lambda_{\bar{\alpha}}%
}\right)  =K\quad\forall\,t\in\lbrack0,T] \label{condizione 2}%
\end{equation}
In order to analyze (\ref{condizione 2}), we need the following Lemma:

\begin{lemma}
\label{lemma 3}Let $a_{m}$ and $b_{0}$ be defined by (\ref{sviluppo 1.1}) and
(\ref{sviluppo 2}) respectively. Under the assumptions of Lemma \ref{lemma 1}
and Lemma \ref{lemma 2}, the function
\[
\arg\left(  4\pi\sqrt{i}\int_{0}^{t}G(t-s)U_{s}\gamma(\underline
{0})\,ds\,e^{-i(T-t)\lambda_{\bar{\alpha}}}\right)
\]
admits the small time expansion%
\begin{gather}
\arg\left(  4\pi\sqrt{i}\int_{0}^{t}G(t-s)U_{s}\gamma(\underline
{0})\,ds\,e^{-i(T-t)\lambda_{\bar{\alpha}}}\right)
=\ \ \ \ \ \ \ \ \ \ \ \ \ \ \ \ \ \ \ \ \ \ \ \ \ \ \ \ \nonumber\\
=\arg\left(  a_{m}\sqrt{i}\right)  +c\sin(\arg b_{0})\,t^{\frac{1}{2}%
}-(T-t)\lambda_{\bar{\alpha}}+o(t^{\frac{1}{2}})\label{sviluppo 3.4}%
\end{gather}
with $c=$ $\frac{B_{m}}{A_{m}}\left\vert \,b_{0}\right\vert $.
\end{lemma}

\begin{proof}
Set $z_{1}=\rho_{1}e^{i\varphi_{1}}$and $z_{2}=\rho_{2}e^{i\varphi_{2}}$; the
first order Taylor expansion of $\arg(z_{1}+z_{2})$ w.r.t. the ratio
$\varepsilon=\frac{\rho_{2}}{\rho_{1}}$ about the point $\varepsilon=0$ is%
\begin{equation}
\arg(z_{1}+z_{2})=\varphi_{1}+\sin(\varphi_{2}-\varphi_{1})\,\varepsilon
+o(\varepsilon) \label{es 2}%
\end{equation}

Using (\ref{carica}) and (\ref{sviluppo 3})we have%
\begin{gather}
\arg\left(  4\pi\sqrt{i}\int_{0}^{t}G(t-s)U_{s}\gamma(\underline
{0})\,ds\,e^{-i(T-t)\lambda_{\bar{\alpha}}}\right)  =\qquad\qquad\qquad
\qquad\qquad\nonumber\\
=\arg\left(  A_{m}a_{m}\,\sqrt{i}\,t^{m+\frac{1}{2}}+B_{m}a_{m}\,\sqrt
{i}\,b_{0}\,t^{m+1}+O(t^{m+\frac{3}{2}})\right)  \,-(T-t)\lambda_{\bar{\alpha
}} \label{sviluppo 3.1}%
\end{gather}
Using (twice) relation (\ref{es 2}), the right hand side of
(\ref{sviluppo 3.1}) can be expanded as
\begin{gather}
\arg\left(  A_{m}a_{m}\,\sqrt{i}\,t^{m+\frac{1}{2}}+B_{m}a_{m}\,\sqrt
{i}\,b_{0}\,t^{m+1}+O(t^{m+\frac{3}{2}})\right)  =\nonumber\\
=\arg\left(  A_{m}a_{m}\sqrt{i}t^{m+\frac{1}{2}}+B_{m}a_{m}\sqrt{i}%
b_{0}t^{m+1}\right)  +O(t)=\nonumber\\
=\arg\left(  a_{m}\sqrt{i}\right)  +c\sin(\arg b_{0})\,t^{\frac{1}{2}%
}+o(t^{\frac{1}{2}}) \label{sviluppo 3.2}%
\end{gather}
Equation (\ref{sviluppo 3.4}) is a straightforward consequence of
(\ref{sviluppo 3.1}) and (\ref{sviluppo 3.2}).
\end{proof}

Lemma \ref{lemma 3} leads us to an asymptotic formulation of the non
controllability condition for small time. From relations (\ref{condizione 2})
and (\ref{sviluppo 3.4}), indeed, we have%
\begin{equation}
K=\arg\left(  a_{m}\sqrt{i}\right)  +\frac{c}{\sqrt{2}}\,t^{\frac{1}{2}%
}-(T-t)\lambda_{\bar{\alpha}}+o(t^{\frac{1}{2}}) \label{condizione 3}%
\end{equation}
where the explicit value $b_{0}=-4\pi\bar{\alpha}\sqrt{i}$ has been taken into
account. Recalling that $c\neq0$, relation (\ref{condizione 3}) is an evident
contradiction we obtained supposing the system (\ref{d-sistema}) to be
non-surjective. This concludes the proof of Theorem \ref{Teorema 1.1} for all
choices of initial states satisfying condition (\ref{3ip}) of Lemma
\ref{lemma 1}.

In the next section we will study an extension of the proof to those cases in
which Lemma \ref{lemma 1} does not applies.

\section{Finite Time Asymptotic for the Charge and Proof of \protect\linebreak
Theorem \ref{Teorema 1.1}}

If condition (\ref{3ip}) does not hold, we may still recover our results by
changing the point in which we perform the expansions of expressions
(\ref{carica}) and (\ref{condizione 1.1}).

To this concern, we consider a radial function $\gamma=\gamma(\left\vert
\underline{x}\right\vert )$ in the space $S(\mathbb{R}^{3})$. From
(\ref{lem1-1}) in Lemma \ref{lemma 1}, it follows that, $U_{t}\gamma
(\underline{0})$ is a $C^{\infty}$-class function represented by%
\begin{equation}
U_{t}\gamma(\underline{0})=\frac{4\pi}{\left(  2\pi\right)  ^{\frac{3}{2}}%
}\int_{0}^{+\infty}k^{2}\mathcal{F}\gamma(k)\,e^{-ik^{2}t}\,dk\label{lem1-1-1}%
\end{equation}
where $\mathcal{F}$ denotes the Fourier transform in $L^{2}(\mathbb{R}^{3})$.
Let us define $f\in L^{2}(-\infty,+\infty)$%
\begin{equation}
f(y)=\left\{
\begin{array}
[c]{l}%
\left(  2\pi\right)  ^{-\frac{1}{2}}\,\,y^{\frac{1}{2}}\,\mathcal{F}%
\gamma(y^{\frac{1}{2}})\quad y\geq0\\
0\qquad y<0
\end{array}
\right.  \label{f}%
\end{equation}
Setting $k^{2}=y$ in the integral (\ref{lem1-1-1}), we can express
$U_{t}\gamma(\underline{0})$ as the Fourier transform of $f$%
\begin{equation}
U_{t}\gamma(\underline{0})=\hat{f}(t)\equiv\int_{0}^{+\infty}f(y)e^{-iyt}%
dy\label{lem1-1-2}%
\end{equation}
Making use of relation (\ref{lem1-1-2}), it is possible to extend $U_{t}%
\gamma(\underline{0})$ to the complex plane as follows%
\begin{equation}
U_{z}\gamma(\underline{0})=\int_{0}^{+\infty}f(y)e^{-iyz}dy;\qquad
z=t+is\label{lem1-1-3}%
\end{equation}
It is well known that, for any $f\in$ $L^{2}(0,+\infty)$, equation
(\ref{lem1-1-3}) defines an holomorphic function in the lower complex half
plane $s<0$ (\cite{Rudin}, Section 19.1). In order to study the limit of
$U_{z}\gamma(\underline{0})$ as $z$ approaches the real axis, we notice that
this function can be expressed as the Fourier transform of the product of
$f(y)$ times $e^{\left\vert y\right\vert s}\in L^{2}(-\infty,+\infty)$%
\[
U_{z}\gamma(\underline{0})=\widehat{\left(  f(y)\cdot e^{\left\vert
y\right\vert s}\right)  }(t)
\]
therefore we have%
\begin{equation}
U_{t+is}\gamma(\underline{0})=\hat{f}\ast\widehat{\left(  e^{\left\vert
\cdot\right\vert s}\right)  }(t)=\frac{1}{\pi}\int_{-\infty}^{+\infty
}U_{t-t^{\prime}}\gamma(\underline{0})\frac{\left\vert s\right\vert
}{\left\vert t^{\prime}\right\vert ^{2}+s^{2}}dt^{\prime}\label{lem1-1-4}%
\end{equation}
with%
\begin{equation}
\frac{1}{\pi}\frac{\left\vert s\right\vert }{\left\vert t^{\prime}\right\vert
^{2}+s^{2}}=\widehat{\left(  e^{\left\vert \cdot\right\vert s}\right)
}(t)\label{lem1-1-5}%
\end{equation}
Taking into account the continuity of $U_{t}\gamma(\underline{0})$, it can be
shown that $U_{z}\gamma(\underline{0})$ is continuous in the set $s\leq0$.

So far we obtained that $U_{t+is}\gamma(\underline{0})$ is holomorphic in
$\left\{  t\in\mathbb{R};\ s<0\right\}  $ and continuous in the closure
$\left\{  t\in\mathbb{R};\ s\leq0\right\}  $. For $\gamma\neq0$, this implies
that the zeroes of function $U_{t}\gamma(\underline{0})$ have to be isolated
points. If $\gamma$ belongs to the domain (\ref{stato iniziale}), the origin
of the real axis is a zero of $U_{t}\gamma(\underline{0})$; therefore a time
$t_{0}>0$, arbitrarily near to the origin, exists such that: $U_{t_{0}}%
\gamma(\underline{0})\neq0$. Proceeding as in Lemma \ref{lemma 1}, it is
possible to obtain for the function $U_{t}\gamma(\underline{0})$ the following
power expansion around $t_{0}$%
\begin{equation}
U_{t}\gamma(\underline{0})=a_{0}+O(\left(  t-t_{0}\right)  );\ a_{0}\neq0
\label{sviluppo 4}%
\end{equation}
with $a_{0}=U_{t_{0}}\gamma(\underline{0})$ and $O(t)\in C^{\infty}%
[0,+\infty)$.

Next, we observe that a simple change of variable%
\[
\tau=t-t_{0}%
\]
and the use of (\ref{soluzione 1}), provide us with an equation for the charge
when the initial time $t=t_{0}$ is assigned%
\begin{equation}
V(0)(\tau)=4\pi\sqrt{i}\int_{0}^{\tau}G(\tau-s)\,U_{t_{0}+s}\gamma
(\underline{0})ds \label{soluzione1.1}%
\end{equation}
Using (\ref{sviluppo 2}) and (\ref{sviluppo 4}), we get the power expansion of
the of the non controllability condition in a right neighbourhood of the point
$t=t_{0}$%
\begin{equation}
K=\arg\left(  a_{0}\sqrt{i}\right)  +\frac{c}{\sqrt{2}}\,\left(
t-t_{0}\right)  ^{\frac{1}{2}}-(T-t)\lambda_{\bar{\alpha}}+o(\left(
t-t_{0}\right)  ^{\frac{1}{2}}) \label{sviluppo 5}%
\end{equation}
with $c\neq0$. As in the previous case, this relation constitutes a
contradiction obtained supposing the system (\ref{d-sistema}) to be non-surjective.

This concludes the proof of Theorem \ref{Teorema 1.1}.

\section{Proof of the Main Result and Final Remarks}

So far, we succeeded in proving that the functional $F(\alpha)$, defined by
(\ref{funzionale})-(\ref{stato iniziale}), belongs to $C^{1}(H_{0}%
^{1}(0,T),\mathbb{C})$ and its derivative, evaluated in $\alpha=0$, is
surjective. The Rank Theorem (see e.g. in \cite{Sontag}, page 336 Theorem 34),
then, applied to our case, implies the existence of a neighborhood of the
origin in $\mathbb{C}$, $I_{0}$, and a $C^{1}$-class map $g:I_{0}\rightarrow
H_{0}^{1}(0,T)$ such that:%
\begin{equation}
F(g(z))=z\quad\forall z\in I_{0} \label{RT}%
\end{equation}
This concludes the proof of Theorem \ref{Teorema 1}.

Our main remark is about the assumptions (\ref{stato iniziale}) on the initial
state. By considering functions with radial symmetry, we are taking into
account only those scattering states which have a null projection along all
spherical harmonics excepting the first one. In this choice there is no loss
of generality. Indeed, those scattering functions whose expansion in spherical
harmonics is:%
\[
\phi(r,\vartheta,\varphi)=\sum_{\substack{l=1\\l\neq0}}^{+\infty}\sum
_{m=-l}^{l}f_{lm}(r)Y_{l}^{m}(\vartheta,\varphi)
\]
exhibit the following characterization:%
\[
U_{t}\phi(\underline{0})=0\ \forall t
\]
Then it follows from definition (\ref{Schroedinger}), and the uniqueness of
solution of the charge equation, that a particle, initially placed in such a
scattering state, $\phi$, and subject to the action of any Hamiltonian of type
$H_{\alpha(t)}$, results to have a null charge and to evolve under the action
of the free propagator; in other words, starting from this initial condition,
the particle doesn't feel the interaction at all. In this case any transfer of
energy is physically impossible. On the other hand our model can be applied to
a realistic situation in which an incoming particle - described by a wave
function of type: $\psi=\phi(r)f(\vartheta,\varphi)$ with $f(\vartheta
,\varphi)$ null outside a cone - is partially trapped into the attractive
potential described by $H_{\bar{\alpha}}$.

In conclusion, we have proved the local controllability of a process of
energy-mass transfer, from scattering to bound states, for a one body quantum
system under the action of a time dependent point interaction.

Further development of this studies may concern the global controllability of
the same process, as well as the inverse problem of finite time ionization.

\section*{Acknowledgments}

I would like to thank Prof. J-M Coron and Prof. R. Figari for useful discussions.

\pagebreak

\appendix

\section*{Appendix A: The Integral Equation}

We recall some basic properties of the integral equations we deal within this
work; some of the relations already used in the previous sections will be here
obtained. A detailed analysis of fractional integral equations, which arise in
the framework of time dependent point interactions in Quantum Mechanics, is
given in \cite{Mik}.

\begin{lemma}
\label{A}Let $\alpha,f\in C(0,T)$; the equation%
\begin{equation}
q(t)+4\sqrt{\pi i}\int_{0}^{t}\frac{\left[  \alpha(s)+\bar{\alpha}\right]
\,\,q(s)}{\sqrt{t-s}}ds=f(t);\qquad\bar{\alpha}\in\mathbb{R}\,
\label{volterra}%
\end{equation}
has an unique continuous solution such that:%
\begin{equation}
\left\Vert q\right\Vert _{C(0,T)}\leq\left\Vert f\right\Vert _{C(0,T)}%
\Gamma(\alpha,T) \label{stima_0}%
\end{equation}
where $\Gamma(\alpha,T)$ is a finite positive constant depending on $\alpha$
and $T$.
\end{lemma}

\begin{proof}
The solution $q(t)$ may be formally expressed using the Picard series:%
\begin{equation}
\left\{
\begin{array}
[c]{l}%
q(t)=\sum_{n=0}^{+\infty}q_{n}(t)\\
q_{0}(t)=f(t)\\
q_{n}(t)=4\sqrt{\pi i}\int_{0}^{t}\frac{\left[  \alpha(s)+\bar{\alpha}\right]
\,q_{n-1}(s)}{\sqrt{t-s}}ds
\end{array}
\right.  \label{Picard}%
\end{equation}
which admit the following estimate:%
\begin{equation}
\sum_{n=0}^{+\infty}\left\vert q_{n}(t)\right\vert \leq\left\Vert f\right\Vert
_{C(0,T)}\left[  1+\sum_{n=1}^{+\infty}\left\vert 4\sqrt{\pi i}\right\vert
^{n}\,\left\Vert \alpha+\bar{\alpha}\right\Vert _{C(0,T)}^{n}\,A_{n}%
\,\pi^{\frac{n}{2}}\,t^{\frac{n}{2}}\right]  \label{stima1}%
\end{equation}
with:%
\[
A_{n}=\left\{
\begin{array}
[c]{l}%
\frac{1}{\left(  \frac{n}{2}\right)  !}\quad n\text{ even}\\
\frac{2^{\left(  \frac{n}{2}+\frac{1}{2}\right)  }}{\Gamma(\frac{1}{2})}%
\frac{1}{n!!}\quad n\text{ odd}%
\end{array}
\right.
\]
being $n!!=1\cdot3\cdot5\cdot\cdot\cdot n$ for $n$ odd. The strong
infinitesimal character of the sequence $A_{n}$ allow this sum to converge
uniformly for any $\alpha\in C(0,T)$ and for any finite time interval $\left[
0,T\right]  $. The existence of a unique solution of (\ref{volterra})
satisfying the estimate (\ref{stima_0}) directly follows from (\ref{stima1}),
with%
\begin{equation}
\Gamma(\alpha,T)=\left[  1+\sum_{n=1}^{+\infty}\left\vert 4\sqrt{\pi
i}\right\vert ^{n}\,\left\Vert \alpha+\bar{\alpha}\right\Vert _{C(0,T)}%
^{n}\,A_{n}\,\pi^{\frac{n}{2}}\,T^{\frac{n}{2}}\right]  \label{Gamma}%
\end{equation}

\end{proof}

In connection with the charge equation in (\ref{Schroedinger}), we claim the
following result:

\begin{lemma}
\label{lemma 6}Let $\psi_{0}\in D(H_{\alpha})$ where $H_{\alpha}$ is the
Hamiltonian associated to a 3-D point interaction placed in the origin and
$\alpha\in\mathbb{R}$; then, the function:%
\begin{equation}
\int_{0}^{t}\frac{U_{s}\psi_{0}(\underline{0})}{\sqrt{t-s}}ds
\label{termine noto 1}%
\end{equation}
is continuous.
\end{lemma}

\begin{proof}
First we consider the case $\alpha<0$.

From definition (\ref{dominio 1}), any function $\psi_{0}\in D(H_{\alpha})$ is
the sum of a regular part plus a bound state term:%
\begin{equation}
\psi_{0}=\varphi+\frac{q}{4\pi\sqrt{2\left\vert \alpha\right\vert }}%
\psi_{\alpha},\,\varphi\in H^{2}(\mathbb{R}^{3}),\,\varphi(\underline
{0})=0,\,q\in\mathbb{C} \label{dominio 1.1}%
\end{equation}
Then the function (\ref{termine noto 1}) disparts in two contributions:%
\begin{equation}
\int_{0}^{t}\frac{U_{s}\psi_{0}(\underline{0})}{\sqrt{t-s}}ds=\int_{0}%
^{t}\frac{U_{s}\varphi(\underline{0})}{\sqrt{t-s}}ds+\frac{q}{4\pi
\sqrt{2\left\vert \alpha\right\vert }}\int_{0}^{t}\frac{U_{s}\psi_{\alpha
}(\underline{0})}{\sqrt{t-s}}ds \label{termine noto 3}%
\end{equation}
We want to prove that (\ref{termine noto 3}) defines a continuous function on
finite time intervals. To this aim we consider the two contributions of
(\ref{termine noto 3}) separately.

The first term of the second member is the one half integral of $U_{t}%
\varphi(\underline{0})$. Using for the state $\varphi$ a representation in
terms of spherical harmonics, holding for all $L^{2}(\mathbb{R}^{3})$
functions, we have:%
\[
\varphi(r,\vartheta,\varphi)=\chi(r)+\sum_{\substack{l=1\\l\neq0}}^{+\infty
}\sum_{m=-l}^{l}f_{lm}(r)Y_{l}^{m}(\vartheta,\varphi)
\]
From the orthogonality relations for $Y_{l}^{m}$, we have%
\[
U_{t}\varphi(\underline{0})=U_{t}\chi(\underline{0})
\]
$\chi\in H^{2}(\mathbb{R}^{3})$ denoting the radial part of $\varphi$. Thus,
the function $U_{t}\varphi$ evaluated in $\underline{x}=0$ may be expressed by
the following Fourier integral%
\begin{equation}
U_{t}\varphi(\underline{0})=\frac{1}{\left(  2\pi\right)  ^{\frac{3}{2}}}%
\int_{\mathbb{R}^{3}}\mathcal{F}\chi(k)\,e^{-i\,k^{2}\,t}\,d\underline{k}
\label{termine noto 3.0}%
\end{equation}
The Fourier transform of the radial function $\chi\in H^{2}(\mathbb{R}^{3})$
has the following characterization: $k^{2}\mathcal{F}\chi(k)\in L^{2}%
(\mathbb{R}^{3})$. Using this relation and the Schwartz inequality, we get an
estimate for the $L^{1}$-norm of $\mathcal{F}\chi$%
\[
\left\Vert \mathcal{F}\chi\right\Vert _{1}=\int_{\mathbb{R}^{3}}%
\frac{\left\vert \mathcal{F}\chi(k)\right\vert \left(  1+k^{2}\right)
}{1+k^{2}}\,d^{3}\underline{k}\leq\left\Vert \left(  1+k^{2}\right)
\mathcal{F}\chi(k)\right\Vert _{2}\left\Vert \frac{1}{1+k^{2}}\right\Vert _{2}%
\]
from which it follows that $\mathcal{F}\chi\in L^{1}(\mathbb{R}^{3})$. This
result guarantees that the function (\ref{termine noto 3.0}), as well as the
first source term $\int_{0}^{t}\frac{U_{s}\varphi(\underline{0})}{\sqrt{t-s}%
}ds$, are continuous.

The second source term of (\ref{termine noto 3}) may be evaluated explicitly
by using the Laplace transform operator $\mathcal{L}$; from definition
(\ref{bound state}), a direct calculation shows that%
\begin{equation}
\mathcal{L}\left[  \int\limits_{0}^{t}\frac{U_{s}\psi_{\alpha}(\underline{0}%
)}{\sqrt{t-s}}\,ds\right]  (p)=\sqrt{\frac{2\pi\left\vert \alpha\right\vert
}{i}}\frac{p^{-\frac{1}{2}}}{p^{\frac{1}{2}}-4\pi\alpha\sqrt{i}}
\label{Ut_bound state2}%
\end{equation}
Here we recall that%
\begin{equation}
\mathcal{L}^{-1}\left[  \frac{p^{-\frac{1}{2}}}{p^{\frac{1}{2}}-4\pi
\alpha\sqrt{i}}\right]  =e^{i\,16\,\pi^{2}\alpha^{2}t}erfc(4\pi\left\vert
\alpha\right\vert \sqrt{i\,t}) \label{termine noto 3.3.0}%
\end{equation}
Then the second contribution of (\ref{termine noto 3}) is%
\begin{equation}
\frac{q}{4\pi\sqrt{2\left\vert \alpha\right\vert }}\int_{0}^{t}\frac{U_{s}%
\psi_{\alpha}(\underline{0})}{\sqrt{t-s}}ds=\frac{q}{4\sqrt{\pi i}%
}e^{i\,16\,\pi^{2}\alpha^{2}t}erfc(4\pi\left\vert \alpha\right\vert
\sqrt{i\,t}) \label{termine noto 3.3.1}%
\end{equation}
which is a continuous function bounded in $\mathbb{R}^{+}$.

The same result holds in the case $\alpha\geq0$, where the only difference
consists in the fact that the operator domain does not include any bound state
of $H_{\alpha}$.
\end{proof}

Relation (\ref{stima_0}) allows us to obtain an estimate for the solution of
charge equation in (\ref{funzionale})%
\begin{equation}
\left\Vert q\right\Vert _{C(0,T)}\leq\left\Vert 4\sqrt{\pi i}\int_{0}^{t}%
\frac{U_{s}\phi(\underline{0})}{\sqrt{t-s}}ds\right\Vert _{C(0,T)}%
\Gamma(\alpha,T) \label{stima 2}%
\end{equation}
where the boundedness of second member is assured by Lemma \ref{lemma 6}.

\begin{quotation}
\textbf{Solving the charge equation}
\end{quotation}

In the particular case $\alpha=0$, the charge $V(0)(t)$ satisfies a fractional
integral equation:%
\begin{equation}
q(t)+4\sqrt{\pi i}\,\bar{\alpha}\int_{0}^{t}\frac{q(s)}{\sqrt{t-s}}%
ds=4\sqrt{\pi i}\int_{0}^{t}\frac{U_{s}\phi(\underline{0})}{\sqrt{t-s}}ds
\label{equazione}%
\end{equation}
whose solution may be explicitly expressed as a functional of the source term;
the Laplace transform of (\ref{equazione}), indeed, gives%
\[
\tilde{q}(p)\left(  1+\frac{4\pi\bar{\alpha}\sqrt{i}}{\sqrt{p}}\right)
=\frac{4\pi\sqrt{i}}{\sqrt{p}}\mathcal{L}\left(  U_{t}\phi(\underline
{0})\right)  (p)\Rightarrow\tilde{q}(p)=\frac{4\pi\sqrt{i}}{\sqrt{p}+4\pi
\bar{\alpha}\sqrt{i}}\mathcal{L}\left(  U_{t}\phi(\underline{0})\right)  (p)
\]
and taking into account the relation:%
\begin{equation}
\mathcal{L}^{-1}\frac{1}{\sqrt{p}+4\pi\bar{\alpha}\sqrt{i}}=\frac{1}{\sqrt{\pi
t}}-(4\pi\bar{\alpha}\sqrt{i})e^{i(4\pi\bar{\alpha})^{2}t}\,erfc(4\pi
\bar{\alpha}\sqrt{it})\equiv G(t) \label{kernel}%
\end{equation}
an explicit expression for the charge is obtained:%
\begin{equation}
V(0)(t)=4\pi\sqrt{i}\int_{0}^{t}G(t-s)\,U_{s}\phi(\underline{0})\,ds
\label{soluzione 1}%
\end{equation}
Due to its similar structure, an analogous expression holds for equation
(\ref{dF0}):%
\begin{equation}
d_{0}V(u)(t)=4\pi\sqrt{i}\int_{0}^{t}G(t-s)\,u(s)\,V(0)(s)\,ds
\label{soluzione 2}%
\end{equation}

When $\phi$ coincide with the bound state $\psi_{\bar{\alpha}}$ of the point
interaction Hamiltonian $H_{\bar{\alpha}}$, the charge equation
(\ref{equazione}) is explicitly solvable. Indeed, the time evolution of a
quantum particle starting in the state $\psi_{\bar{\alpha}}$ and moving under
the action of $H_{\bar{\alpha}}$, is%
\[
\psi(t)=e^{-itH_{\bar{\alpha}}}\psi_{\bar{\alpha}}=e^{-it\lambda_{\bar{\alpha
}}}\psi_{\bar{\alpha}}%
\]
Comparing this expression with (\ref{dominio 1}), we get the explicit form of
the charge%
\begin{equation}
q(t)=4\pi\sqrt{2\left\vert \bar{\alpha}\right\vert }\,e^{-it\lambda
_{\bar{\alpha}}} \label{soluzione 3}%
\end{equation}
Relations (\ref{soluzione 1}) and (\ref{soluzione 3}) imply%
\begin{equation}
\sqrt{i}\int_{0}^{t}G(t-s)\,U_{s}\psi_{\bar{\alpha}}(\underline{0}%
)\,ds=\sqrt{2\left\vert \bar{\alpha}\right\vert }\,e^{-it\lambda_{\bar{\alpha
}}} \label{soluzione 4}%
\end{equation}
Replacing $t$ with $T-t$, the previous formula can be written as%
\[
\sqrt{i}\int_{0}^{T-t}G(T-t-s)\,U_{s}\psi_{\bar{\alpha}}(\underline
{0})\,ds=\sqrt{2\left\vert \bar{\alpha}\right\vert }\,e^{-i\left(  T-t\right)
\lambda_{\bar{\alpha}}}%
\]
Moreover, setting $s^{\prime}=T-s$ in the integral at the l.h.s. and dividing
by $\sqrt{i}$, we get%
\begin{equation}
\int_{t}^{T}G(s^{\prime}-t)U_{T-s^{\prime}}\psi_{\bar{\alpha}}(\underline
{0})\,ds^{\prime}=\sqrt{\frac{2\left\vert \bar{\alpha}\right\vert }{i}%
}e^{-i(T-t)\lambda_{\bar{\alpha}}} \label{carica inversa 2}%
\end{equation}
This relation has been used in order to obtain condition (\ref{condizione 1}).

In order to justify the use of Laplace transform in deriving relations
(\ref{soluzione 1}) and (\ref{soluzione 2}), our next task is to prove the following

\begin{lemma}
\label{lemma 0} Let $q(t)$ be the solution of the charge equation:%
\[
q(t)+4\sqrt{\pi i}\int_{0}^{t}\frac{\left[  \alpha(s)+\bar{\alpha}\right]
\,\,q(s)}{\sqrt{t-s}}ds=4\sqrt{\pi i}\int_{0}^{t}\frac{U_{s}\psi
_{0}(\underline{0})}{\sqrt{t-s}}ds
\]
with $\alpha\in L^{\infty}(\mathbb{R})$, $\psi_{0}\in D(H_{\alpha(0)})$. Then,
the Laplace transform $\mathcal{L}q(p)$ exists and is analytic at least in the
open half plane of $\mathbb{C}$ defined by the condition:%
\begin{equation}
p\in\mathbb{C}:\operatorname{Re}p>16\pi^{2}\left\Vert \alpha\right\Vert
_{\infty}^{2} \label{dominio laplace}%
\end{equation}

\end{lemma}

\begin{proof}
Consider the function: $q^{\prime}(t)=e^{-p\,t}q(t)$ with $p\in\mathbb{C}$ and
$\operatorname{Re}p>0$; it satisfies the equation:%
\begin{equation}
q^{\prime}(t)+4\sqrt{\pi i}\int_{0}^{t}\frac{\alpha(s)\,e^{-p(t-s)}%
\,q(s)}{\sqrt{t-s}}ds=4\sqrt{\pi i}\int_{0}^{t}\frac{e^{-p(t-s)}}{\sqrt{t-s}%
}f(s)\,ds \label{carica 2.1}%
\end{equation}
where the function $f$ in the nonhomogeneous term is given by:%
\begin{equation}
f(t)=e^{-p\,t}U_{t}\,\psi_{0}(\underline{0}) \label{termine noto 3.3}%
\end{equation}
Using relations (\ref{dominio 1.1}) and (\ref{termine noto 3.0}) of \ the
previous proof, we may write (\ref{termine noto 3.3}) in the following form:%
\begin{equation}
f(t)=e^{-p\,t}\left[  \frac{1}{\left(  2\pi\right)  ^{\frac{3}{2}}}%
\int_{\mathbb{R}^{3}}\mathcal{F}\varphi_{r}(k)\,e^{-i\,k^{2}\,t}%
\,d\underline{k}+\frac{q}{4\pi\sqrt{2\left\vert \alpha\right\vert }}U_{t}%
\psi_{\alpha}(\underline{0})\right]  \label{termine noto 3.4}%
\end{equation}
with $\mathcal{F}\varphi_{r}(k)\in L^{1}(\mathbb{R}^{3})$; again an explicit
calculation shows that:%
\begin{multline*}
U_{t}\psi_{\alpha}(\underline{0})=\mathcal{L}^{-1}\left[  \sqrt{\frac
{2\left\vert \alpha\right\vert }{i}}\frac{1}{\left(  p^{\frac{1}{2}}%
+4\pi\left\vert \alpha\right\vert \sqrt{i}\right)  }\right]  =\\
=\sqrt{\frac{2\left\vert \alpha\right\vert }{i}}\left[  \frac{1}{\sqrt{\pi t}%
}-(4\pi\left\vert \alpha\right\vert \sqrt{i})e^{i(4\pi\alpha)^{2}t}%
\,erfc(4\pi\left\vert \alpha\right\vert \sqrt{it})\right]
\end{multline*}
from which we deduce that $f(t)\in L^{1}(0,+\infty)$.

Next we apply the Young's inequality:%
\[
\left\Vert f\ast g\right\Vert _{1}\leq\left\Vert f\right\Vert _{1}\left\Vert
g\right\Vert _{1}%
\]
for convolutions of $L^{1}(0,+\infty)$ functions, to the equation
(\ref{carica 2.1}), and obtain the following estimate:%
\[
\left\Vert q^{\prime}\right\Vert _{1}\left(  1-4\sqrt{\pi}\left\Vert
\alpha\right\Vert _{\infty}\int_{0}^{+\infty}\left\vert \frac{e^{-p\,t}}%
{\sqrt{t}}\right\vert dt\right)  \leq4\sqrt{\pi}\left\Vert f\right\Vert
_{1}\int_{0}^{+\infty}\left\vert \frac{e^{-p\,t}}{\sqrt{t}}\right\vert dt
\]
which provides an effective bound for the norm $\left\Vert q^{\prime
}\right\Vert _{1}$ if the coefficient $\left(  1-4\sqrt{\pi}\left\Vert
\alpha\right\Vert _{\infty}\int_{0}^{+\infty}\left\vert \frac{e^{-p\,t}}%
{\sqrt{t}}\right\vert dt\right)  $ is positive. Recalling that, for
$\operatorname{Re}p>0$, holds the equality:%
\[
\int_{0}^{+\infty}\left\vert \frac{e^{-p\,t}}{\sqrt{t}}\right\vert
dt=\sqrt{\frac{\pi}{\operatorname{Re}p}}%
\]
we get the condition:%
\begin{equation}
\left(  1-4\sqrt{\pi}\left\Vert \alpha\right\Vert _{\infty}\int_{0}^{+\infty
}\left\vert \frac{e^{-p\,t}}{\sqrt{t}}\right\vert dt\right)  >0\Rightarrow
1>4\sqrt{\pi}\sqrt{\frac{\pi}{\operatorname{Re}p}}\left\Vert \alpha\right\Vert
_{\infty}\Rightarrow\operatorname{Re}p>16\pi^{2}\left\Vert \alpha\right\Vert
_{\infty}^{2} \label{4ip}%
\end{equation}
Following the same line, it's easy to prove that, if the condition (\ref{4ip})
holds, the partial derivatives of the function $q^{\prime}$ w.r.t the real and
the imaginary part of $p$ - both given by:%
\[
-te^{-p\,t}q(t)
\]
are bounded by integrable functions of $t$:%
\[
\int_{0}^{+\infty}\left\vert te^{-p\,t}q(t)\right\vert dt\leq\int_{0}%
^{+\infty}\left\vert te^{-\left(  16\pi^{2}\left\Vert \alpha\right\Vert
_{\infty}^{2}+\varepsilon\right)  \,t}q(t)\right\vert dt<\infty
\]
Then $e^{-p\,t}q(t)$ is $C^{1}$ integrable w.r.t. $t\in\left[  0,+\infty
\right)  $ for any $p$ in the domain (\ref{dominio laplace}); moreover, in the
same hypothesis, its partial derivatives w.r.t. $p$ are bounded by measurable
functions of $t$. This allows us to conclude that the Laplace integral:%
\[
\mathcal{L}q(p)=\int_{0}^{+\infty}q(t)e^{-p\,t}\,dt
\]
defines a $C^{1}$ class function for $p$ in the domain (\ref{dominio laplace}).
\end{proof}

\end{document}